\documentstyle [12pt]{article}

\begin {document}
\title{Rigid $\aleph_\epsilon$-saturated models of superstable
theories\footnote{Partially supported by the United States Israeli
Science Foundation Pub nu.569}} 
\author {Ziv Shami and Saharon Shelah}

\newtheorem {theorem}{Theorem}
\newtheorem {lemma}{Lemma}      
\newtheorem {definition}{Definition}
\newtheorem {proof_of_lemma}{Proof of Lemma}
\newtheorem {proof_of_theorem}{Proof of Theorem}
\newtheorem {fact}{Fact}
\newtheorem {conclusion}{Conclusion}

\newtheorem {remark}{Remark}

\newtheorem {claim}{Claim}
\newtheorem {claim_proof}{Proof of Claim}
\newtheorem {example}{Example}
\newtheorem {notations}{Notations}

\maketitle

\begin{abstract}
 In a countable superstable NDOP theory, the existence of a rigid $\aleph_{\epsilon}$-saturated model implies the existence of $2^\lambda$ rigid $\aleph_{\epsilon}$-saturated models of power $\lambda$ for every $\lambda > 2^{\aleph_0}$.
\end{abstract}

\newbox\noforkbox\newdimen\forklinewidth
\newbox\forkbox\newdimen\forklinewidth
\newbox\unionbox\newdimen\forklinewidth

\forklinewidth=0.3pt   

\setbox0\hbox{$\textstyle\cup$}
\setbox2\hbox{$\textstyle\uplus$}
\setbox3\hbox{$\textstyle\bigcup$}

\setbox1\hbox to \wd0{\hfil\vrule width \forklinewidth depth \dp0
                        height \ht0 \hfil}
\wd1=0 cm

\setbox\noforkbox\hbox{\box1\box0\relax}
\setbox\forkbox\hbox{\box1\box2\relax}
\setbox\unionbox\hbox{\box1\box3\relax}

\def\unionstick{\mathop{\copy\noforkbox}\limits}
\def\unionstcka{\mathop{\copy\forkbox}\limits}

\def\nonfork#1#2_#3{#1\unionstick_{\textstyle #3}#2}
\def\fork#1#2_#3{#1\unionstcka_{\textstyle #3}#2}

\def\nonforkin#1#2_#3^#4{#1\unionstick_{\textstyle #3}^{\textstyle #4}#2}     
\def\card #1 {\vert #1 \vert}


\section{Introduction}
         
Ehrenfuecht conjectured that given a theory $T$, the class of cardinals 
for which $T$ has a rigid model is quite well behaved.
Shelah refuted Ehrenfeucht's conjecture showing that this class can be quite 
complicated.\\
\\
In this paper we deal with problems related to this question in the context
of stability. More specifically, we will study the existence of stable rigid 
models which satisfy an additional saturation property (note that if no 
\\saturation property is required, then a very simple example of a stable rigid
model can be found, namely the 
model whose language is $\{P_n \vert n<\omega \}$ and consists of the disjoint
union of the $P_n$-s, each of which has exactly one element.)\\

We will give a partial solution to the following questions:\\
1) What classes of superstable theories have a rigid $\aleph_\epsilon$-saturated 
model?\\
2) Assuming that there exits a rigid $\aleph_\epsilon$-saturated model, what 
can be said about the number of $\aleph_\epsilon$-saturated 
models, or perhaps even about the number of rigid $\aleph_\epsilon$-saturated models, 
in large enough cardinality?\\
In section 3 we consider two properties of a superstable theory $T$:\\
1) $T$ is strongly deep.\\
2) $T$ does not admit a nontrivial nonorthogonal automorphism of some saturated
model.\\
We prove that 1) is a necessary condition for the existence of a rigid
$\aleph_\epsilon$-saturated model, and that 2) is a sufficient condition.\\
In section 4, we assume that $T$ is a superstable NDOP theory. We prove that
2) is actually equivalent to the existence of a rigid 
$\aleph_\epsilon$-saturated model. We can then conclude
that the existence of a single rigid $\aleph_\epsilon$-saturated model implies
that $T$ has $2^\lambda$ such models for every sufficiently large $\lambda$. 
 
In this paper, the notations will be very similar to Shelah's notations in [Sh-C]. 
$T$ will denote any complete stable theory with no finite models in some 
language $L$. $\kappa, \lambda, \mu$ will denote cardinals.\\
We work in some huge saturated model $\cal M$. Sets $A,B,C$,... will be
subsets of $\cal M$, with cardinality strictly less than the 
cardinality of $\cal M$. $\bar a, \bar b, \bar c,...$ will denote finite tuples
of $\cal M$.$M$, $N$ will always be elementary submodels of $\cal M$.
$p,q,r$ will denote types, usually complete over some set A.
$S(A)$ will denote the set of complete types over $A$. Also, for a tuple $\bar b$ and a set $A$, 
$\bar b/A$ denote the type of $\bar b$ over $A$.     
 
\section{Building a dimensionally diverse $\aleph_{\epsilon}$-saturated model}
In this section, T denotes any stable theory. We give a brief outline of some
standard constructions following [Sh-C].   
$\kappa(T)$ herein will denote the smallest infinite cardinal $\kappa$ 
such that there is no chain $\{p_\alpha\in S(A_\alpha):\alpha <\kappa\}$ 
such that for all $\alpha <\beta <\kappa$, $p_\beta$ is a forking extension
of $p_\alpha$. 
Recall that $M$ is $F^a_\mu$-saturated if for every $A\subseteq M$ such that $\card {A} < \mu$, 
every strong type over $A$ is realized in $M$. $M$ is $F^a_\mu$-prime
over $A$ if $M$ is $F^a_\mu$-saturated and for every $F^a_\mu$-saturated model $N$ such that $A\subseteq N$
there is an elementary embedding of M into N over A. 
We say that $M$ is $\aleph_\epsilon$-saturated ($\aleph_\epsilon$-prime over A)
if $M$ is $F^a_\omega$-saturated ($F^a_\omega$-prime over A).

\begin{definition}\em
 We say that an $\aleph_\epsilon$-saturated model M (of a superstable T) is 
 $dimensionally$ $diverse$, if for any two stationary regular types
 $p$, $q$ over finite subsets of M, $p \perp q$ if and only if      
 $dim(p,M)\neq dim(q,M)$.
\end{definition}
We also recall the following standard definition:

\begin{definition}\em
 For superstable T, we say that T is $multidimensional$ if for every cardinal $\alpha$, 
 there are nonalgebraic $p_{i},$ $i<\alpha$ which are pairwise orthogonal. 
\end{definition}
Our aim in this section is to prove the following standard theorem.
\begin{theorem}
 Let T be superstable, and let $\mu,\delta$ be cardinals such that $\aleph_\delta=\delta,$ $\mu<\delta$ and $2^{\card{T} }<\delta$,
 then\\  
 1) There exists a $\mu$-saturated model M of cardinality $\delta$ 
 which is dimensionally diverse  
(in particular if $\mu \geq \aleph_{1}$,
  M is $\aleph_{\epsilon}$-saturated).\\ 
2) If T is multidimensional, then for every increasing sequence of cardinals 
 $\bar{\mu}=\langle\mu_{\alpha}\vert \alpha < \delta\rangle,$ 
 $\mu_{\alpha}\in(\mu,\delta)$ there exists a $\mu$-saturated and 
 $\aleph_{\epsilon}$-saturated model M of cardinality $\delta$ 
 which is dimensionally diverse, and such that for every $\alpha < \delta$ there
 is a stationary and regular(=s.r.) type $p_{\alpha}$ (over a finite subset of M,) 
 such that for every s.r. type $\bar{p}$,
 $dim(\bar{p},M)=\mu_{\alpha}$ if and only if $\bar{p}\not\perp p_{\alpha}$.    
 Moreover, for every s.r. type $\bar{p}$ (over a finite subset of M) 
 $dim(\bar{p},M)=\mu_{\beta}$ for some $\beta<\delta$
       
\end{theorem}
\begin{fact}
 Let $M=\bigcup_{i<\alpha} {M_{i}}$,
 where $\langle M_{i} \vert\  i<\alpha \rangle$
 is an increasing and continuous sequence of $\kappa$-saturated models, let
 $A\subseteq M_{0}$ such that $\card{A} < \kappa$, and let $p\in S(A)$ be a stationary
 regular type. Then $dim(p,M)=dim(p,M_{0})+\sum_{i<\alpha}\  
 dim(p_{i},M_{i+1})$, where 
 $p_{i}$ is the stationarization of p to $M_{i}$.     
\end{fact}

\begin{fact}
 If $cf(\lambda)\ge \kappa(T), M$ is $F^a_\lambda$-prime over $A$, 
 and $I\subseteq M$ is an indiscernible sequence over $A$, then 
 $\card I \le \lambda$. 
\end{fact}

Following are two claims which we will use in our proof of Theorem 1 (although
weaker versions thereof would suffice.)
\begin{claim}
 Suppose $cf(\mu)\geq\kappa=\kappa(T)$ and $M\models T\ is\ F^{a}_\kappa$-saturated. Let $A\subseteq$ $\cal M$ 
 such that $\card{A} \leq\mu$, and suppose $M^{+}$ is $F_{\mu}^{a}$-prime over $M\cup A$.
 In addition, let $B\subseteq M^{+}$ such that $\card{B} < \kappa$, and suppose $p\in S(B) $ is stationary 
 and $\lambda=dim(p,M^{+}) > \mu$. Then $p\not\perp M$. 
 Moreover, if p is regular then there is a stationary regular type 
 $q\in S(B^{*}),$ where $B^{*}\subseteq M$ and 
 $\card {B^{*}} < \kappa $ such that $dim(q,M)=\lambda$ and $q\not\perp p$.  
\end{claim}

\begin{claim_proof}\em
Left to the reader.
\end{claim_proof}

\begin{claim}
 Suppose $cf(\mu)\geq\kappa=\kappa(T)$ and $M\models T$ is $F_{\kappa}^{a}-saturated$. Let $p_{i}\in S(M)$, $i < \alpha$ be pairwise
 orthogonal and let $E=\bigcup_{i < \alpha}{E_{i}}$ where $E_{i}$ is a Morley sequence of $p_{i}$.
 In addition, suppose N is $F_{\mu}^{a}$-prime over
 $M\cup E$, and let $B\subseteq N$ such that $\card{B} < \kappa$. If $q\in S(B)$ is stationary 
 and regular and $dim(q,N) >\mu$, then 
 $q\not\perp M$.   
\end{claim}

\begin{claim_proof}\em
 Assume, for a contradiction, that $q\perp M$. Let $M^+$ be 
 $F_\mu^a$-prime over $M\cup B$ with $M^+ \prec N$, and let $\tilde q \in S(M^+)$
 be the stationarization of $q$. Then either $dim(q,M^+) > \mu$ or $dim(\tilde{q},N) > \mu$, by Fact 1.
 Now, there exists $S\subseteq \alpha$ with $\card{S} < \kappa$, such that 
 $p_i\perp tp(M^+/M)$ for all $i \in \alpha \setminus S$. Thus $dim(q,M^+) > \mu$ 
 contradicts Claim 1, and $dim(\tilde{q},N) > \mu$ contradicts the above and
 Fact 2.   
\end{claim_proof}

\begin{fact}
 If $cf(\delta)\ge \kappa(T)$, and $\langle M_i \vert i<\delta \rangle$
 is an elementary chain of $\lambda$-saturated models, then
 $M^*=\bigcup_{i<\delta}{M_i}$ is $\lambda$-saturated.
 In particular, if T is superstable, then the union of any elementary chain
of $\lambda$-saturated models is $\lambda$-saturated.  
\end{fact}

Theorem 1 easily follows from the above:

\begin{proof_of_theorem}\em
 Sketch of proof: 
 Define by induction an increasing elementary chain $\langle M_{\alpha}\vert \alpha < \delta\rangle$
 of $\mu$-saturated models. At the $\alpha^{th}$ step, choose a non-algebraic 
 stationary type $p_\alpha$ (over a finite set) which is orthogonal to 
 $p_\beta$ for all $\beta<\alpha$, and define $M_\alpha$ to be $\mu$-prime
 over $\bigcup_{\beta<\alpha} M_\beta \cup Dom(p_\alpha) \cup I_\alpha$,
 where $I_\alpha$ is a Morley sequence of $p_\alpha$ with cardinality 
 $\mu_\alpha$. Now let $M^*=$ $\bigcup_{\alpha<\delta} M_\alpha$. From Fact 1,
 Claim 1 and Claim 2 it follows that $M^*$ realizes the desired dimensions 
 ($\mu_\alpha$), and by Fact 3, $M^*$ is $\mu$-saturated.    
\end{proof_of_theorem}

\section{The existence of a rigid $\aleph_{\epsilon} $-saturated model.} 
 In this section T is assumed to be superstable. We say that a superstable
 theory is $strongly$ $deep$ if the depth of every type is positive (if and only if the 
 depth of every type is infinity). We prove that whenever T has a rigid 
 $\aleph_\epsilon$-saturated model, T is strongly deep. We also introduce 
 the notion of a $nontrivial$ $nonorthogonal$ $automorphism$ and prove
 that if some saturated model does not have such an automorphism, then
 in arbitrarily large cardinality, T has a maximal number of rigid $\aleph_{\epsilon}$-saturated 
 models. (i.e. $2^\lambda$ such models in cardinality $\lambda$.)

\begin{definition}\em
We say that T is $strongly$ $deep$ if for every $\aleph_{\epsilon}$-saturated 
model M, and for every (without loss of generality regular) type $p\in S(M)$, if N is $\aleph_{\epsilon}$
-prime over M$\cup\{a \}$ where a$\models p$, then $q\perp M$ for some $q\in S(N)$.     
\end{definition}

\begin{lemma}
 Let $N_{0}$ be $\aleph_{\epsilon}$-saturated, let $p\in S(N_{0})$ be regular,
 and let $\langle \bar{e}_{i}\vert i<\alpha \rangle$ be a Morley sequence of p. 
 Suppose $ N^{+}$ is $\aleph_{\epsilon}$-prime over $N_{0}\cup\bigcup_{i<\alpha}
 {\bar{e}_{i}}$. Then the following are equivalent:\\\\
 (i) There is $p_{1}\in S(N_{1})$ such that $p_{1}\perp N_{0}$,
 where $N_{1}$ is $\aleph_{\epsilon}$-prime over 
   $N_{0}\cup\bar{e_{0}}$.\\
 (ii) There is $p^{+}\in S(N^{+})$ such that
 $p^{+}\perp N_{0}$.      

\end{lemma}

\begin{proof_of_lemma}\em
 (i)$\rightarrow$(ii) We may assume that $N_{1}\prec N^{+}$. If there is such $p_{1}$, choose $p^+$ to be the nonforking
 extension of $p_{1}$ to $N^{+}$.\\
(ii)$\rightarrow$(i) Let us assume that  
 such a $p^{+}$ is given, then $p^{+}$ is strongly based on some finite subset A of
 $N^+$. Therefore
 $tp(A/N_0\cup \bigcup_{i<\alpha}{\bar{e}_i } )$ is $F_{\aleph_{0}}^{a}$-isolated, 
 and we may assume $\alpha=n< \omega$. Now, if every $p\in S(N_{1})$ is nonorthogonal
 to $N_{0}$, then by induction on n we have $p^{+}\not\perp N_{0}$
 (recall that the depth of parallel types is the same) which is a contradiction.    
\end{proof_of_lemma}

\begin{theorem}
 If T has an $\aleph_{\epsilon}$-saturated model N such that $\card{Aut(N)} <2^{\aleph_{0}} $,
 then T is strongly deep.
\end{theorem}

\begin{proof_of_theorem}\em
 Suppose not, then there is some depth 0 type p. Let $N_{0}\prec N$ 
 be $\aleph_{\epsilon}$-prime over $\emptyset$ (without loss of generality $p\in S(N_{0})$ ). 
 Let $I\subseteq N$ be a maximal Morley sequence of p (without loss of generality $I$ is infinite,)  
 and let M be a maximal model such that $N_{0}\prec M\prec N$ and $\nonfork{M}{I}_{N_{0}}$ (so M is $\aleph_{\epsilon}$-saturated).
 We claim that N is $\aleph_{\epsilon}$-minimal over $M\cup I$.   
 Otherwise, let $P\prec N$ be $\aleph_{\epsilon}$-prime over $M\cup I$, so for
 some $\bar{b}\in N$, $\bar{b}/P $ is a non-algebraic regular type.

 $\underline{Claim\  1}\ \ $ $\bar{b}/P \perp N_{0}\ $. 
 
 Indeed, if not then (by [Sh-C])
 there is $\bar{c}\in N$ such that $\bar{c}/P$ is regular, $\bar{c}/P \not\perp \bar{b}/P$ 
 and $\nonfork{\bar{c}}{P}_{N_{0}}$. Thus $\{\bar{c},M,I\}$ is independent over $N_{0}$, 
in contradiction to the definition of M.

 $\underline{Claim\  2}$ For all $\bar{d}\in N$ if $\bar{d}/M \perp N_{0}$
 then $\bar{d}\in M$.

 Indeed, let $\bar{d} $ satisfy the above, so $\nonfork{\bar{d}}{I}_{M}$ and  
 $\nonfork{M}{I}_{N_{0}}$. Therefore $\nonfork{M\cup\bar{d}}{I}_{N_{0}}$,   
   so by the maximality of M $\ \bar{d}\in M$. 

 $\underline{Claim\ 3}\ \ \bar{b}/P \perp M$.

 Otherwise $\bar{b}/P \not\perp M$, so there is $\bar{c}\in N $ with $\bar{c}/P$       
 regular such that $\nonfork{\bar{c}}{P}_{M}$ and $\bar{c}/P\not\perp\bar{b}/P$. 
 Hence by claim 1 $\bar{c}/P \perp N_{0}$,  
 so $\bar{c}/M \perp N_{0}$. But then applying Claim 2, we get $\bar{c}\in M$, which is
 a contradiction.

  Now, continuing to prove the theorem, we recall that $tp(I/M)$ does not fork over $N_{0}$, so Claim 3 and lemma 1 imply that
 for all $\bar{e}\in I$, Depth $(\bar{e}/M) > 0$. But 
 $Depth(\bar{e}/M)=Depth(\bar{e}/N_{0})=0$, a contradiction.    
 So, we have shown that $N$ is $\aleph_{\epsilon}$-minimal over $M\cup I$, 
 and therefore $\aleph_{\epsilon}$-prime over $M\cup I$.
 By the uniqueness of $\aleph_{\epsilon}$-prime models
 and the fact that $\nonfork{M}{I}_{N_{0}}$, we conclude that every 
 permutation of $I$ induces an automorphism of $N$, thus 
 $\card{Aut(N)} \geq 2^{\aleph_{0}}$, which is a contradiction.       
\end{proof_of_theorem}

\begin{conclusion}
 If T has an $\aleph_{\epsilon}$-saturated model N such that 
 $\card{Aut(N)} < 2^{\aleph_{0}}$, then T is multidimensional. 
 (because even $Depth(T)>0$ implies that T is multidimensional.)       
\end{conclusion}

\begin{definition}\em
 We say that $\sigma\in Aut(M) $ ($M$ is $\aleph_{\epsilon} $-saturated) is a 
 $nontrivial$ $nonorthogonal$ $automorphism$ (=n.n.a.) if for any nonalgebraic $p\in S(M)$,
 $p\not\perp\sigma(p)$, and $\sigma\neq id$.
\end{definition}

\begin{remark}\em
 $\sigma\in Aut(M)$ ($M$ as above) is a n.n.a. if and only if the unique extension of it
 $\sigma^{eq}$ to $M^{eq}$ is a n.n.a. 
\end{remark}

\begin{theorem}
If $\delta=\aleph_{\delta}>\beta\geq 2^{\card{T} }$ and the saturated model of 
cardinality $\beta$ does not have a nontrivial nonorthogonal automorphism,  
then T has $2^{\delta}$ rigid $\beta$-saturated models of cardinality
$\delta$. 
\end{theorem}

\begin{proof_of_theorem}\em
 It is enough to show that every dimensionally diverse $\beta$-saturated 
 model of cardinality $\delta$ is rigid. This is indeed enough, because 
using Theorem 1 we may then take a dimensionally diverse $\beta$-saturated 
model $N$ of cardinality $\delta$. So $N$ is rigid and by Conclusion 1 T is 
multidimensional. Then, by Theorem 1, for every subset $D$ of 
$(\beta,\delta)$ which consists of cardinals, we may choose a 
$\beta$-saturated model which realizes exactly the dimensions in $D$ (in the
 sense of Theorem 1, part 2,) from which the theorem follows.\\   
 So let M be a dimensionally diverse $\beta$-saturated model of 
cardinality $\delta$. And let $\sigma\in Aut(M)$. Assume by way of
contradiction that $\sigma\neq id$. We define by induction an increasing 
chain of elementary submodels of M: Let $M_{0}$ be a saturated model of 
cardinality $\beta$, such that $\sigma\vert M_{0}\neq id$. For all $n<\omega$
let $M_{i+1}$ be $F_{\beta}^{a}$-prime over $\bigcup_{i\in Z}{\sigma^{i}(M_{n})}$.
Now, according to Fact 3, $M_{n}$ is the saturated model of cardinality $\beta$ 
or all $n<\omega$. And if $M_{\omega}=\bigcup_{i<\omega}{M_{i}}$, then 
$M_{\omega}$ is also the saturated model of cardinality $\beta$. Clearly,
 $\sigma\vert M_{\omega}$ is a nontrivial automorphism of $M_{\omega}$, 
so by the assumption of the theorem, there is a 
 regular type $p\in S(M_{\omega})$ such that $p\perp \sigma(p)$. Hence for
 some finite $A\subset M_{\omega}$, p is strongly based on A,
 so $dim(p\vert A,M)=dim(\sigma(p\vert A),M)$, which contradicts the fact 
 that M is dimensionally diverse.

\end{proof_of_theorem}

\section{A characterization for superstable NDOP countable theories.}
 In this section T is assumed to be a supertstable NDOP countable theory.
 We will use [SH-401] in order to get the following characterization: such a theory
 has a rigid $\aleph_\epsilon$-saturated model if and only if every $\aleph_\epsilon$
 -staurated model does not have a nontrivial nonorthogonal automorphism.    
 Then we will use this characterization to show that the existence of a single
 rigid $\aleph_\epsilon$-saturated model implies the existence of a maximal
 number of such models in every sufficiently large cardinality.\\
 We work in $\cal M$$^{eq}$. 

\begin{definition}\em
 We say that $A$ is $almost$ $finite$ if $A$ is contained in the algebraic closure
 of some finite set. 
\end{definition}

 \subsection {The $L_{\infty,\aleph_{\epsilon}}(d.q)$-Characterization Theorem 
  of [SH-401]}
 
 $\ $
 
 Let $M_{0},M_{1}$ be $\aleph_{\epsilon}$-saturated. We say that they are
 $L_{\infty,\aleph_{\epsilon}}(d.q)-equivalent$, if there is a family $\cal F$ 
which satisfies the following:\\

$\ $
1) Each $f\in$ $\cal F$ is an elementary partial map from $M_{0}$ to $M_{1}$
such that Dom(f) is almost finite. 
 
 $\ $

 2) $\cal F $ is closed under restriction.
 
 $\ $ 
 
 3) For every $f\in\cal F$, and every $\bar{a_{l}}\in M_{l}$  ($\bar{a_{l}} $ are finite sequences,
    $l=0,1)$ there exists $g\in$ $\cal F$ such that $f\subseteq g$,
 $acl(\bar{a_{0}})\subseteq Dom(g)$, and $acl(\bar{a_{1}})\subseteq Rang(g)$.    
 
 $\ $

 4) Whenever $tp(e/A)$ (where A is a almost finite,) is stationary and regular, then for 
some almost finite $A^{*}\supseteq A$,
    if $f\cup\{\langle e_{0},e_{1}\rangle\}\in$ $\cal F$
    and $tp(e_{0}/Dom(f))$ is conjugate to the stationarization 
    of $tp(e/A)$ to $A^{*}$, then
  
    $dim(p,Dom(f),I_{0})=dim(f(p),Rang(f),I_{1})$, where p denotes $e_{0}/Dom(f)$,
    
    $I_{0}=\{e\in M_{0}\vert\ f\cup \{\langle e,e_{1}\rangle\}\in$ ${\cal F}$ $\}$, and   
    
    $I_{1}=\{e\in M_{1}\vert\ f\cup \{\langle e_{0},e\rangle\}\in$ ${\cal F}$ $\}$. \\ \\
Conditions 4) is the main assumption, roughly speaking, it implies that: whenever ${\cal F}$ sends
a stationary and regular type p to a type q, then the full structure of dimensions above p will be the
same as the full structure of dimensions above q.

We can now state the Characterization Theorem for $\aleph_\epsilon$-saturated
 models $M_{0},$ $M_{1}$: \\ \\
$M_{0}, M_{1}$ $are$ $isomorphic$ $if$ $and$ $only$ $if$ $they$ $are$ $L_{\infty,\aleph_{\epsilon}}(d.q)$-
 $equivalent.$

\begin{theorem}
 If T has a nontrivial nonorthogonal automorphism 
 of some $\aleph_{\epsilon}$
 -saturated model, then every $\aleph_{\epsilon}$-saturated model of T is not
 rigid. 
\end{theorem} 

\begin{proof_of_theorem}\em
 By Remark 1, without loss of generality $T=T^{eq}$.  
 Suppose $N_{0}$ is an $\aleph_{\epsilon}$-saturated model which has a n.n.a. 
 $\sigma_{0}$, and let $a_{0}\neq a_{1}$ be in $N_{0}$ such that $\sigma_{0}(a_{0})=a_{1}$. 
 Let $N$ be some $\aleph_{\epsilon}$-saturated model (of $T$), we will show that
 $N$ is not rigid. Choose $b_{0},b_{1}\in N$ such that $tp(b_{0},b_{1})=tp(a_{0},a_{1})$.\\
 Let $\cal F$ be the family of all partial elementary maps $f$ from $(N,b_{0})$ to $(N,b_{1})$ 
 with an almost finite domain, such that for some partial elementary map $\tau$,  
 with almost finite domain, $({\sigma_0\vert {A_0}})\tau=\tau\sigma$ $\}$,
 where $A_0=rang(\tau)$. By [SH-401] it is sufficient to show that $\cal F$ satisfies 1) - 4) above. 

 1), 2) are immediate. 3) is also immediate by the fact that $N$ is 
 $\aleph_{\epsilon}$-saturated. 
 To show 4), we will prove the following claim: if $tp(e/A)$ is stationary and regular, 
 where $A$ is almost finite, then it can be replaced by a nonforking extension
 of it, denoted again by $tp(e/A)$ (where $A$ is almost finite and $A\cup\{e\}\subseteq N$,) 
 such that if $\sigma_{i}(A)=A^{'}$ and  $\sigma_{i}(e)=e_{i}$ for 
 $\sigma_i\in$ $\cal F$ ($i=0,1)$ then $\fork{e_{0}}{e_{1}}_{A^{'}}$. 4) will follow from
 this, as it implies that for all $f\in$ $\cal F$ with $Dom(f)=A$, if $e^*_0=e,\ e^*_1=f(e), M_0=M_1=N$ 
 and $I_0, I_1$
 are defined as in 4) of subsection 4.1 (with $e^*_i$ instead of $e_i$,) then     
 $dim(tp(e^*_0/A),A,I_0)=dim(tp(e^*_1/f(A)),f(A),I_1)=1$.
 
 Proceeding to prove this claim, we note that by Theorem 2, T is without loss 
 of generality strongly deep. Let $M_{0}$ be an 
 $\aleph_{\epsilon}$-prime model over $\emptyset$, such that $tp(e/M_{0})$
 is a nonforking extension of $tp(e/A)$, let
 $M_{0}^{+}$ be $\aleph_{\epsilon}$-prime over $M_{0}\cup \{e\}$ and let $e^{+}$     
 be such that $tp(e^{+}/M_{0}^{+})\perp M_{0}$, with $tp(e^{+}/M_{0}^{+})$ non 
 algebraic and regular. Now let $B\subseteq M_{0}^{+}$ be finite such that 
 $tp(e^{+}/M_{0}^{+})$ is strongly based on $B$, so there is some finite 
 $C\subseteq M_{0}$ such that $tp(B/acl(C\cup \{e\}))\vdash 
 tp(B/M_{0}\cup \{e\})$. We may assume without loss of generality that 
 $M_{0}^{+}\subseteq N$ (because the $\aleph_{\epsilon}$-prime model over 
 the union of $M_{0}$ and a countable set is also $\aleph_{\epsilon}$-prime 
 over $\emptyset$ (see [Sh-C]). 
 
 Now, by the way condition 4) was stated we may also assume that $A\supseteq C$.  
 Assume, towards a contradiction, that $\nonfork{e_{0}}{e_{1}}_{A^{'}}$, and
 let $\sigma_{i}^{+}$ be in $\cal F$ with $\sigma_i \subseteq \sigma_{i}^{+}$
 and $Dom(\sigma_{i}^{+})\supseteq
 B\cup acl(A\cup \{e\})$, for i=0,1. Choose $M^{'}_0$ to be an $\aleph_{\epsilon}$-prime model over $\emptyset$, 
 such that $A'\subseteq M^{'}_0$ and $\nonfork{\{e_{0},e_{1}\}}{M_{0}^{'}}_{A'}$, and set $B_{i}=\sigma_{i}(B)$.
 Then by the fact that $tp(e_{0},e_{1}/M_{0}^{'})$ determines the type 
 $tp(acl(\{e_{0} \} \cup A^{'}), acl(\{ e_{1}\} \cup A^{'})/ M_{0}^{'})$, and the fact
 that $tp(B/acl(\{e\}\cup A))\vdash tp(B/\{e\}\cup M_{0})$, we conclude that
 we may choose $M_{0}^{*}$ and $M_{1}^{*}$ which are  $\aleph_{\epsilon}$-prime 
 over $M_{0}^{'}\cup \{e_{0}\}$ and $M_{1}^{'}\cup \{e_{1}\}$ respectively, 
 and such that $B_{i}\subseteq
 M_{i}^{*}$ for i=0,1. So if we set $q=tp(e^{+}/B)$ and $q_{i}=\sigma_{i}^{+}(q)$
 for $i=0,1$, then by the definition of $\cal F$, $q\not\perp q_{i}$ for i=0,1. 
 Let $\bar{q_{i} }\in S(M_{i}^{*})$ be a nonforking extension of $q_{i}$. 
 Since $\bar{q_{i} }\perp M_{0}^{'}$ (as $M_{i}^{*}$ is $\aleph_{\epsilon}$-prime 
 over $M_{i}^{'}\cup 
 \{e_{i}\}\cup B_{i}$ and unique upto isomorphism over $M_{i}^{'}\cup 
 \{e_{i}\}\cup B_{i}$,) and since $\nonfork{M_{0}^{*}}{M_{1}^{*}}_{M_{0}^{'}}$, 
 we conclude that $q_{0}\perp q_{1}$. Thus $q$, $q_{0}$ and $q_{1}$ contradict the 
 fact that nonorthogonality is an equivalence relation on stationary regular types. 
\end{proof_of_theorem}

\begin{theorem}
If T has a rigid $\aleph_{\epsilon}$-saturated model then
for every cardinal $\lambda\geq {(2^{\aleph_{0}})}^{+}$ T has $2^{\lambda}$
rigid $\aleph_{\epsilon}$-saturated models of power $\lambda$. 
\end{theorem}

Before proving the theorem, we will need a combinatorial lemma. To that end,
we first introduce the following notations.  

\begin{notations}\em
1)  $\cal T$ denotes a subtree of $^{<\omega}\lambda$. Let $\eta, \nu\in$ $\cal T$, then\\
(i) $\nu^-=\eta$ if and only if $\nu$ is a successor of $\eta$.\\
(ii) $lg(\eta)$ denotes the length of $\eta$.\\
(iii) ${\cal T}_\nu = \{ \eta \in {\cal T} | \nu \triangleleft 
\eta \}$, ${\cal T}^+_\nu = \{ \eta \in {\cal T} | \nu \triangleleft
\eta$ or $\nu=\eta\}$.

2)  $\cal R$ $=\langle ( N_{\eta}, a_\eta)| \eta \in \cal T \rangle$ denotes an 
    $\aleph_\epsilon$-representation (see Chapter X, Def. 5.2 in [Sh-C]).\\
3) For an $\aleph_\epsilon$-representation $\cal R$=$\langle ( N_{\eta}, a_\eta)| \eta \in \cal T \rangle$ 
   $E_{ \cal R }$ denotes the equivalence relation on $\cal T$ defined by:
$E_{ \cal R }(\eta, \eta')$ if and only if there exists $v \in {\cal T}$ such that $\eta 
= v {\scriptstyle \wedge} \langle \alpha \rangle, \eta' = v {\scriptstyle \wedge} 
\langle \alpha ' \rangle$ for some 
$\alpha, \alpha^{'} $, and $tp(a_\eta/N_\nu)=tp(a_{\eta^{'}}/N_\nu)$.\\

\end{notations}

\begin{definition}
\begin{enumerate}\em
\item[{\rm (1)}] 
We say that a subtree $\cal T$ of $^{<\omega}\lambda$ is $\mu-wide$ if 
for all \\ $\eta\in$ $\cal T$,  $\vert\{\nu\in \cal{T} \vert \nu^-=\eta\}\vert\ge\mu$.  
  
\item[{\rm (2)}] Suppose ${\cal T}_0, {\cal T}_1 \subseteq {}^{< \omega} 
\lambda$ 
are subtrees, we say that ${\cal T}_0, {\cal T}_1$ are $\mu-equivalent$ 
if there exists $\bar{A_i}=\langle A^\rho_i \vert \rho\in {\cal T}_i \rangle\ (i=0,1) $
such that\\  
(a) for i = 0, 1 and all $\rho\in {\cal T}_i$, $\vert A^\rho_i \vert < \mu$.\\ 
(b) for i = 0, 1 $A^\rho_i\subseteq {\cal T}_\rho$, and \\

(c) ${\cal T}_0 \setminus\bigcup_{\tau\in {U_0}}{({\cal T}_0)_\tau}$, ${\cal T}_1 \setminus\bigcup_{\tau\in {U_1}}{({\cal T}_1)_\tau}$  
are isomorphic as partial orders, where 
$U_i=\bigcup_{\rho\in {\cal T}_i}{A_i^\rho}$. 

\item[{\rm (3)}] A tree ${\cal T} \subseteq {}^{< \omega} \lambda$ is 
called 
$\mu-strongly$ $rigid$ if for every $\eta$ and 
$\alpha_0 < \alpha_1 < \lambda \quad$ (such that $\eta{\scriptstyle \wedge}\langle{\alpha_i}\rangle\in{\cal T}$) 
${\cal T}^+_{\eta^{\scriptstyle \wedge} \langle \alpha_0 \rangle}, {\cal T}^+_{\eta^{\scriptstyle \wedge} \langle \alpha_1 \rangle}$ 
are not $\mu$-equivalent, and ${\cal T}$ is $\mu$-wide.

\item[{\rm (4)}] We say that an $\aleph_\epsilon$-representation 
$\cal R$ $=\langle ( N_{\eta}, 
a_\eta)| \eta \in \cal T \rangle$ is 

\item[{\rm (i)}]  $of$ $maximal$ $width$, if for every $\nu \in \cal T$, 
$\{ tp (a_\eta / N_{\eta^-} )| \ \eta\in {\cal T} , \eta^- = \nu \}$ is a maximal set 
of pairwise orthogonal regular types (modulo equality) such that 
$tp (a_\eta /  
N_{\eta^-}) \bot N_{\nu^-}$. 

\item[{\rm (ii)}] $equally$ $divided$, if $| \eta \ {\scriptstyle \wedge}  \langle \alpha \rangle / E_{\cal R} | = 
| \eta \ {\scriptstyle \wedge} \langle \beta \rangle / E_{\cal R}| $ for all 
$\eta$, $\alpha$, and $\beta$ such that $\eta{\scriptstyle \wedge}\langle{\alpha}\rangle, \eta{\scriptstyle \wedge}\langle{\beta}\rangle\in {\cal T} $.

\item[{\rm (iii)}] $\mu-wide$, 
if $\vert\eta/E_{\cal R}\vert\ge \mu$ for all $\eta\in$ $\cal T $.  
\end{enumerate}
 
\end{definition}
Before stating the following lemma, we would like to point out that much stronger
versions of it have been proved by Shelah.
\begin{lemma}
 For every $\lambda > \mu$ there exist $2^\lambda$ trees ${\cal T}$ 
$\subseteq \lambda^{<\omega}$ which are $\mu-$ strongly rigid and $\mu$-
 nonequivalent, of power $\lambda$.
\end{lemma}

\begin{proof_of_lemma}\em 
 First, it is enough to show the existence of a single such tree, in which
 the root has $\lambda$ many successors. So, let us construct such a tree:\\

 a) Let $\langle A_n \vert n<\omega \rangle$ be a partition of the set of positive
 natural numbers, such that for all $n<\omega\  A_n$ is infinite and $min (A_n) > n$.
 Also, suppose $A_n=\{ k^n_l \vert l<\omega \}$, where $k^n_i<k^n_j$ for $i<j$.\\
 b) Let $ h: \lambda \rightarrow \lambda $ be surjective, such that 
 for all $\alpha < \lambda$,  $\card{ h^{-1}(\alpha) } =\lambda$.\\
 c) Let ${\cal T }_0= ^{<\omega} \lambda=\{\eta_i \vert i<\lambda \}$ 
    ($\eta_i\neq\eta_j$ for $i\neq j$).\\       
 d) For every ordinal $i$, let $t_i=\{\nu \vert \nu$ is a strictly decreasing sequence
    of ordinals $ <\omega+i\}$.\\
 e) Now, let us define our tree: ${\cal T }^*=\{\rho\in ^{<\omega} \lambda\ \vert$
    for all $k < lg(\rho)$, if $(*)(k,\rho)$ then $\rho(k) < \mu \}$, where
 $(*)(k,\rho)$ is the following statement:\\\\  If $n<\omega$ is the unique natural number such that $k\in A_n$,
 and $l(*),i(*)$ are such that $k=k^n_{l(*)}$ and $\rho\vert n=\eta_{i(*)}$, then 
 $\langle h(\rho(k^n_l)) \vert l<l(*) \rangle \not\in t_{i(*)}$.\\\\
 f) The following ranks are useful. Let $(n,S,\eta)$ be a triple, with $n<\omega$, $S\subseteq ^{<\omega}\lambda$ is a subtree, and $\eta\in S$.  
 We define an ordinal rank $rk^n[\eta,S]$, by:\\
 i) $rk^n[\eta,S] \geq 0$ for every $\eta\in S$.\\
 ii) $rk^n[\eta, S]\geq \alpha+1$ if  
     there exist $\{ \eta^*_i \vert i<\mu^+ \}$, where $\eta^*_i$  
     are distinct elements of $S$, each extending $\eta$ and satisfying $lg(\eta^*_i)\in A_n$, and $rk^n[\eta^*_i,S]\geq\alpha$.\\
 g) It can be easily verified that:\\
  i) If $S\subseteq T\subseteq ^{<\omega}\lambda$ are $\mu$- equivalent subtrees, 
     $\eta\in S$, and $n<\omega$, then $rk^n[\eta,S]=rk^n[\eta,T]$.\\  
  ii) If $i<\lambda$, $n<\omega$ and $lg(\eta_i)=n$, where $\eta_i\in {\cal T^*}$
      then $rk^n[\eta_i,$ ${\cal T^*}]=\omega+i$.\\
 $\ $ 
 Now, by g) we conclude that $\cal T^*$ is $\mu$-strongly rigid.    
\end{proof_of_lemma}

{\bf Fact 5(Chapter X, [Sh-C])} suppose that $T$ is a superstable NDOP theory. 
Let ${\cal R}_i=$  $\langle (N^i_{\eta}, a^i_\eta) | \eta \in 
{\cal T}_i \rangle$ be $\aleph_\epsilon$-representations which are $\aleph_1$-wide, 
let $M_i$ be $\aleph_\epsilon$-prime over 
$\bigcup_{\eta\in {\cal T}_i}{N^i_\eta}$
and let $\sigma: M_0 \to M_1$ be 
an isomorphism. Then there are ${\cal T}^\ast_i \subseteq {\cal T}_i$ 
such that ${\cal T}^\ast_i, {\cal T}_i$ 
are $\aleph_1$ -equivalent, for $i = 0, 1$ and an isomorphism 
$\tilde \sigma : {\cal T}^\ast_0 
\to {\cal T}^\ast_1$ (of partial orders) such that

for all $\eta_i\in {\cal T}^*_i\ (i=0,1)$ if $p^i_{\eta_i}=tp(a^i_{\eta_i}/
N^i_{\eta^-_i})$, then
$ \sigma (p^0_{\eta_0}) {\not\perp} p^1_{\eta_1}$ implies 
$\tilde \sigma ((\eta_0 / 
E_{{\cal R}_0})\cap {\cal T}^*_0) = (\eta_1 / E_{{\cal R}_1})\cap {\cal T}^*_1.$

\begin{proof_of_theorem}\em
By Lemma 2, for every $\lambda > 2^{\aleph_0}$ there 
exist $2^\lambda$ trees ${\cal T} \subseteq {}^{< \omega} \lambda$ which 
are $2^{\aleph_0}$-strongly rigid and $2^{\aleph_0}$-non-equivalent, of power 
$\lambda$.  

We will show that for every such ${\cal T}$, if ${\cal R}=\langle (M_{\eta},e_\eta) | \eta \in {\cal T} \rangle$ 
is an equally divided $\aleph_\varepsilon$-representation of maximal width, and $M$ is $\aleph_\varepsilon$-prime 
over it, then $M$ is rigid.

This is sufficient because by Fact 5, any two $\aleph_\varepsilon$-prime models over $2^{\aleph_0}$-non-equivalent 
trees are non-isomorphic.

So, assume $\sigma \in {\rm Aut} (M)$, we must show $\sigma = {\rm id}$: otherwise 
by Theorem 4, there exists a non-algebraic regular type $p^\ast \in S (M)$ such that $p^\ast \bot \sigma p^\ast$. Since $T$ 
has NDOP and $M$ is $\aleph_\varepsilon$-minimal over $\bigcup_{\eta\in {\cal T}}{M_\eta}$,
we have $p^\ast {\not\perp}
M_{\eta_0}$ for some $\eta_0\in{\cal T}$. By the fact that $M$ is $\aleph_\varepsilon$-prime 
over an $\aleph_\varepsilon$-representation of maximal width, $p^\ast {\not\perp} p_{{\eta_{0}}^\wedge {\langle\alpha^\ast \rangle}}$ 
for some $\alpha^\ast$. Therefore 
$\sigma p_{\eta^\wedge_0 \langle \alpha^\ast\rangle}$ 
${\not\perp} \sigma p^\ast \bot p^\ast \not\perp 
p_{\eta^\wedge_0 \langle \alpha^\ast\rangle}$. 
In particular, $(\ast) \quad \sigma p_{\eta^{\wedge \langle \alpha^\ast\rangle}_0} 
\bot p_{\eta^{\wedge \langle \alpha^\ast\rangle}_0}$. But by Fact 5, there are ${\cal T}^\ast_0, {\cal T}^\ast_1 
\subseteq {\cal T}$ (without loss of generality
$\eta_0 {\scriptstyle\wedge} \langle \alpha^\ast \rangle \in {\cal T}^\ast_0$), and an isomorphism 
of partial 
orders $\tilde \sigma: {\cal T}^\ast_0 \to {\cal T}^\ast_1$, such that 
${\cal T}^\ast_i (i = 0, 1)$ are 
$\aleph_1$-equivalent to ${\cal T}$, and $\eta^\ast$ such that 
$\sigma p_{\eta_0 \wedge \langle \alpha^\ast \rangle} 
\not\perp p_{\eta^\ast}$. So by $(\ast)$ $p_{\eta^\ast} \neq 
p_{\eta^{\wedge 
\langle \alpha^\ast\rangle}_0}$ and  
by Fact 5, $\tilde \sigma (\eta {\scriptstyle\wedge} \langle\alpha^\ast\rangle /E_{\cal R} ) = 
\eta^\ast / E_{\cal R} \neq (\eta_0 {\scriptstyle\wedge} \langle\alpha^\ast\rangle) / E_{\cal R}$, 
in contradiction to the fact 
that ${\cal T}$ is even $2^{\aleph_0}$-strongly rigid (here we used the fact that ${\cal T}$ is equally divided).

\end{proof_of_theorem}

\begin{example}\em
 Let $L=\{E,f\}$ and let $M$ be the following $L$-structure:\\
 $\vert M \vert=Z\times\omega$, $E$ is the equivalence relation defined by
 $E[(m_0,k_0),(m_1,k_1)]$ if and only if $m_0=m_1$, and $f:\vert M \vert \rightarrow 
 \vert M \vert$ is any function with the following properties:\\
 i) For all $\ (m,k)\in\vert M \vert\ \ f(m,k)=(m+1,k^{'}) $ for some $k^{'} $.\\
 ii) For all $\ (m,k)\in\vert M \vert\ \ f^{-1}(m,k) $ is infinite.\\
 It is not hard to see that $T=Th(M)$ is an $\omega$-stable NDOP theory      
 in which any saturated model does not have a n.n.a., and therefore $T$ has
 $2^\lambda$ rigid $\aleph_\epsilon$-saturated models of cardinality $\lambda$
 for every $\lambda > 2^{\aleph_0}$.  
\end{example}

Ziv Shami, Department of Mathematics, University of Notre Dame, Notre Dame, IN 46556-5683, USA\\
E-mail address: Ziv.Shami.1@nd.edu\\
\\ 
Saharon Shelah, Institute of Mathematics, The Hebrew University, Jerusalem, Israel\\
E-mail address: Shelah@math.huji.ac.il

\end{document}